\documentclass{gtart_h}

\def\ifplaintex{\expandafter\ifx\csname documentclass\endcsname\relax}

\ifplaintex 
\else
\fi

\def\gt{{\mathsurround=0pt\it $\cal G\mskip-2mu$eometry \&\ 
$\cal T\!\!$opology}}        %  journal title in recommended style

\def\gtm{{\mathsurround=0pt\it $\cal G\mskip-2mu$eometry \&\ 
$\cal T\!\!$opology $\cal M\mskip-1mu$onographs}}    %  for monographs

\def\gtp{{\mathsurround=0pt\it $\cal G\mskip-2mu$eometry \&\ 
$\cal T\!\!$opology $\cal P\!$ublications}}  % GT publications

\def\recd{{\small Received:\qua\receiveddate\ifx\reviseddate\relax
\else\qquad Revised:\qua\reviseddate\fi\par}} 

\def\volumenumber#1{\def\thevolumenumber{#1}}
\def\volumeyear#1{\def\thevolumeyear{#1}}
\def\volumename#1{\def\thevolumename{#1}}
\def\papernumber#1{\def\thepapernumber{#1}}
\def\pagenumbers#1#2{\def\startpage{#1}\def\finishpage{#2}}
\def\published#1{\def\publishdate{#1}}
\def\received#1{\def\receiveddate{#1}}
\def\revised#1{\def\reviseddate{#1}}
\def\accepted#1{\def\accepteddate{#1}}
\def\asciititle#1{\def\theasciititle{#1}}
\def\covertitle#1{\def\thecovertitle{#1}}
\def\coverauthors#1{\def\thecoverauthors{#1}}
\def\asciiauthors#1{\def\theasciiauthors{#1}}

\def\coverauthors#1{\def\thecoverauthors{#1}}
\long\def\asciiabstract#1{\long\def\theasciiabstract{#1}}
\def\asciikeywords#1{\def\theasciikeywords{#1}}

\let\\\par
\let\thevolumenumber\relax\let\thepapernumber\relax
\let\thevolumeyear\relax\let\startpage\relax
\let\finishpage\relax\let\publishdate\relax\let\receiveddate\relax
\let\reviseddate\relax\let\accepteddate\relax\let\theasciititle\relax
\let\thecovertitle\relax\let\theasciiauthors\relax
\let\theasciiabstract\relax\let\theasciikeywords\relax
\let\thecoverauthors\relax
\let\thecoverauthors\relax\let\theerratum\relax\let\theasciiemail\relax
\let\theshortauthors\relax\let\theshorttitle\relax

\def\startpage{1}\def\finishpage{15}\def\thepapernumber{77}

\volumenumber{2}
\volumename{Proceedings of the Kirbyfest}
\volumeyear{1999}

\long\def\maketitlep{   % start of definition of \maketitlep

\count0=\startpage

\gtm\nl        %   GT mongraphs (top left) 
{\small Volume \thevolumenumber: \thevolumename\nl 
\ifx\theerratum\relax\else Erratum \erratumnumber\nl\fi
Pages \startpage--\finishpage\nl}

\vglue 0.1truein   % top margin

{\parskip=0pt\leftskip 0pt plus 1fil\def\\{\par\smallskip}{\ifplaintex\large
\else\Large\fi\bf\thetitle}\par\medskip}   
\vglue 0.05truein 

{\parskip=0pt\leftskip 0pt plus 1fil\def\\{\par}{\sc\theauthors}
\par\medskip}%
 
\vglue 0.03truein 

{\small\leftskip 25pt\rightskip 25pt{\bf Abstract}\stdspace\theabstract

{\bf AMS Classification}\stdspace\theprimaryclass
\ifx\thesecondaryclass\relax\else; \thesecondaryclass\fi\par
{\bf Keywords}\stdspace \thekeywords\par}\vglue 7pt

}   % end of definition of \maketitlep

\font\phead=cmsl9 scaled 950
\font=cmsl9 scaled 1050
\font\pnum=cmbx10 scaled 913
\font\lnum=cmbx10 
\font\pfoot=cmsl9 scaled 950
\font=cmsl9 scaled 1050
\ifplaintex
\headline{\vbox to 0pt{\vskip -4.5mm\line{\small\phead\ifnum
\count0=\startpage ISSN 1464-8997 (on line)
1464-8989 (printed) \hfill {\pnum\folio}\else\ifodd\count0\def\\{ }% 
\ifx\theshorttitle\relax\thetitle\else\theshorttitle\fi\hfill{\pnum\folio}
\else\def\\{ and }{\pnum\folio}\hfill\ifx\theshortauthors\relax\theauthors
\else\theshortauthors\fi\fi\fi}\vss}}
\footline{\vbox to 0pt{\vglue 0mm\line{\small\pfoot\ifnum\count0=\startpage
Published \publishdate:\qua\copyright\ \gtp\hfill\else
\gtm, Volume \thevolumenumber\ (\thevolumeyear)\hfill\fi}\vss
}}
\else
\makeatletter
\def\@oddhead{{\small\lhead\ifnum\count0=\startpage ISSN 1464-8997 (on line)
1464-8989 (printed) \hfill {\lnum\number\count0}\else\ifodd\count0
\def\\{ }\ifx\theshorttitle\relax \thetitle \else\theshorttitle\fi\hfill
{\lnum\number\count0}\else\def\\{ and }{\lnum\number\count0}
\hfill\ifx\theshortauthors\relax 
\theauthors\else\theshortauthors\fi\fi\fi}}\def\@evenhead{@oddhead}
\def\@oddfoot{\small\lfoot\ifnum\count0=\startpage Published \publishdate:\qua\copyright\ \gtp\hfill\else
\gtm, Volume \thevolumenumber\ (\thevolumeyear)\hfill\fi}
\def\@evenfoot{@oddfoot}
\makeatother
\fi

\let\maketitlepage\maketitlep
\let\makeshorttitle\maketitlepage
\let\maketitle\maketitlepage

\newwrite\gtoutfile
\long\gdef\makeheadfile{  %%% start of definition of \makeheadfile
{\def\\{, }\def\s{ }
\immediate\openout\gtoutfile head.xxx
\immediate\write\gtoutfile{Proxy-for: \ifx\theasciiauthors\relax
\theauthors\else\theasciiauthors\fi\s<\ifx\theasciiemail\relax\theemail\else\theasciiemail\fi>}
\immediate\write\gtoutfile{\noexpand\\}
\immediate\write\gtoutfile{Authors: \ifx\theasciiauthors\relax
\theauthors\else\theasciiauthors\fi}
{\def\\{ }\immediate\write\gtoutfile{Title: \ifx\theasciititle\relax
\thetitle\else\theasciititle\fi}}
\immediate\write\gtoutfile{Subj-class: GT or SG, GR etc}
\immediate\write\gtoutfile{MSC-class: \theprimaryclass\ifx\thesecondaryclass\relax\else, \thesecondaryclass\fi}
\immediate\write\gtoutfile{Journal-ref: Geom. Topol. Monogr. \thevolumenumber\s
(\thevolumeyear) \startpage-\finishpage}
\immediate\write\gtoutfile{Comments: Published by Geometry and Topology Monographs at}
\immediate\write\gtoutfile{\s\s\s  http://www.maths.warwick.ac.uk/gt/GTMon\thevolumenumber/paper\thepapernumber.abs.html}
\immediate\write\gtoutfile{\noexpand\\}
\immediate\write\gtoutfile{}
\ifx\theasciiabstract\relax
\immediate\write\gtoutfile{\theabstract}\else
\immediate\write\gtoutfile{\theasciiabstract}\fi
\immediate\write\gtoutfile{}
\immediate\write\gtoutfile{\noexpand\\}
\immediate\write\gtoutfile{}
\immediate\closeout\gtoutfile}}  %%% end of definition of \makeheadfile

\def\maketitlepage{\maketitlep\makeheadfile}
\let\makeshorttitle\maketitlepage
\let\maketitle\maketitlepage

\volumenumber{7}
\volumename{Proceedings of the Casson Fest}
\volumeyear{2004}
\papernumber{1}
\pagenumbers{1}{26}
\received{22 October 2003}
\revised{29 March 2004}
\accepted{10 March 2004}
\published{17 September 2004}

\usepackage{amsmath,amssymb,latexsym}

\def \ol{\overline}
\def \pd {\partial}
\def \tl {\widetilde}
\def \wtl{\widetilde}
\def \what{\widehat}
\def \lraw {\longrightarrow}
\def \daw {\downarrow}

\def \ga {\alpha}
\def \gb {\beta}

\def \gd {\delta}

\def \gs {\sigma}
\def \cC {{\mathcal C}}

\def \BC {{\mathbb C}}
\def \BH {{\mathbb H}}
\def \BQ {{\mathbb Q}}
\def \BR {{\mathbb R}}
\def \BZ {{\mathbb Z}}

\def \fG {\mathfrak{G}}
\def \Coker {{\rm Coker\,}}
\def \Ker {{\rm Ker\,}}
\def \Hom {{\rm Hom\,}}
\def \la {\langle}
\def \ra {\rangle}
\def \sm {-} %\setminus

\def \beq {\begin{equation}}
\def \eeq {\end{equation}}

\newcommand{\cd}{\operatorname{cd}}

\newtheorem{theorem}{Theorem}[section]
\newtheorem{corollary}[theorem]{Corollary}
\newtheorem{conjecture}{Hope}[section]
\newtheorem{proposition}[theorem]{Proposition}
\newtheorem{addendum}[theorem]{Addendum}
\newtheorem{lemma}[theorem]{Lemma}

\begin{document}
\title{Poincar\'e duality in dimension 3}
\covertitle{Poincar\noexpand\'e duality in dimension 3}
\asciititle{Poincare duality in dimension 3}
\author{C\,T\,C Wall}
\coverauthors{C\noexpand\thinspace T\noexpand\thinspace
C Wall}
\asciiauthors{CTC Wall}
\begin{abstract}\footnotesize
The paper gives a review of progress towards extending the Thurston
programme to the Poincar\'e duality case. In the first section, we fix
notation and terminology for Poincar\'e complexes $X$ (with fundamental
group $G$) and pairs, and discuss finiteness conditions.\nl
For the case where there is no boundary, $\pi_2$ is non-zero if and
only if $G$ has at least 2 ends: here one would expect $X$ to split
as a connected sum. In fact, Crisp has shown that either $G$ is a free
product, in which case Turaev has shown that $X$ indeed splits, or $G$ is
virtually free. However very recently Hillman has constructed a Poincar\'e
complex with fundamental group the free product of two dihedral groups
of order 6, amalgamated along a subgroup of order 2.\nl
In general it is convenient to separate the problem of making the
boundary incompressible from that of splitting boundary-incompressible
complexes. In the case of manifolds, cutting along a properly embedded
disc is equivalent to attaching a handle along its boundary and then
splitting along a 2--sphere. Thus if an analogue of the Loop Theorem
is known (which at present seems to be the case only if either $G$ is
torsion-free or the boundary is already incompressible) we can attach
handles to make the boundary incompressible. A very recent result of
Bleile extends Turaev's arguments to the boundary-incompressible case,
and leads to the result that if also $G$ is a free product, $X$ splits
as a connected sum.\nl
The case of irreducible objects with incompressible boundary can be
formulated in purely group theoretic terms; here we can use the recently
established JSJ type decompositions. In the case of empty boundary the
conclusion in the Poincar\'e duality case is closely analogous to that
for manifolds; there seems no reason to expect that the general case
will be significantly different.\nl
Finally we discuss geometrising the pieces. Satisfactory results follow
from the JSJ theorems except in the atoroidal, acylindrical case, where
there are a number of interesting papers but the results are still far
from conclusive.\nl
The latter two sections are adapted from the final chapter of my survey
article on group splittings.
\end{abstract}

\asciiabstract{%
The paper gives a review of progress towards extending the Thurston
programme to the Poincare duality case. In the first section, we fix
notation and terminology for Poincare complexes X (with fundamental
group G) and pairs, and discuss finiteness conditions.  For the case
where there is no boundary, pi_2 is non-zero if and only if G has at
least 2 ends: here one would expect X to split as a connected sum. In
fact, Crisp has shown that either G is a free product, in which case
Turaev has shown that X indeed splits, or G is virtually free. However
very recently Hillman has constructed a Poincare complex with
fundamental group the free product of two dihedral groups of order 6,
amalgamated along a subgroup of order 2.  In general it is convenient
to separate the problem of making the boundary incompressible from
that of splitting boundary-incompressible complexes. In the case of
manifolds, cutting along a properly embedded disc is equivalent to
attaching a handle along its boundary and then splitting along a
2-sphere. Thus if an analogue of the Loop Theorem is known (which at
present seems to be the case only if either G is torsion-free or the
boundary is already incompressible) we can attach handles to make the
boundary incompressible. A very recent result of Bleile extends
Turaev's arguments to the boundary-incompressible case, and leads to
the result that if also G is a free product, X splits as a connected
sum.  The case of irreducible objects with incompressible boundary can
be formulated in purely group theoretic terms; here we can use the
recently established JSJ type decompositions. In the case of empty
boundary the conclusion in the Poincare duality case is closely
analogous to that for manifolds; there seems no reason to expect that
the general case will be significantly different.  Finally we discuss
geometrising the pieces. Satisfactory results follow from the JSJ
theorems except in the atoroidal, acylindrical case, where there are a
number of interesting papers but the results are still far from
conclusive.  The latter two sections are adapted from the final
chapter of my survey article on group splittings.}

\asciikeywords{Poincare complex, splitting, loop theorem, incompressible,
  JSJ theorem, geometrisation}
\address{Department of Mathematical Sciences, University of 
Liverpool\\Liverpool, L69 3BX, UK}
\email{ctcw@liv.ac.uk}
\maketitle
\cl{\small\it Dedicated to Andrew J Casson on the occasion of his
60th birthday}

\section{Preliminaries}\label{secprelim}

In \cite{PC1} I defined Poincar\'e complexes as CW--complexes satisfying
the strongest (global) form of the Poincar\'e duality theorem that holds
for manifolds.

Suppose given a connected complex $Y$, with fundamental group
$G:=\pi_1(Y)$ having group ring $R:=\BZ G$ (we fix these notations
throughout this article); and a homomorphism $w\co G\to\{\pm 1\}$. The
ring $R=\BZ G$ admits the involutory anti-automorphism $(\sum a_g
g)^*=\sum a_g g^{-1}$. This allows us to regard any left $R$--module
as a right module and vice-versa; by default, we use `$R$--module' for
right $R$--module. Write $C_*(Y)$ for the chain complex of the universal
cover $\tl{Y}$, regarded as a complex of free $R$--modules. Then for any
(right) $R$--module $B$, we set $H^*(Y;B):=H^*(\Hom_R(C_*(Y),B))$ and
${}^tH_*(Y;B):=H_*(C_*(Y)\otimes_R B)$. Here the affix $t$ is to emphasise
that we used the homomorphism $w$ to transfer the given right module
structure on $B$ to a left module structure.

If we are given a class $[Y]\in {}^tH_n(Y;\BZ)$ such that, for all
$r\in\BZ$, cap product with $[Y]$ induces an isomorphism
$$[Y]\frown \co H^r(Y;R) \lraw {}^tH_{n-r}(Y;R\otimes\BZ^t);$$
then we call $Y$ a connected $\mathit{PD}^n$ complex with fundamental class $[Y]$.
According to \cite[Lemma 1.1]{PC1}, it follows that $[Y]$ is unique (up
to sign) and that for any $r\in \BZ$ and any $R$--module $B$, we have an
isomorphism $[Y]\frown\co H^r(Y;B)\to {}^tH_{n-r}(Y;B)$. We say that $Y$
is orientable if $w$ is trivial; $[Y]$ then defines an orientation. We
will usually assume orientability.

The above definition contains no explicit finiteness condition. However,
following \cite[page 222]{brownbk}, we may argue as follows. The
homology and cohomology groups are those of a complex $C_*$, say,
of $R$--modules. The functors $B\to {}^tH_k(Y;B)$ commute with
direct limits. Since $[Y]\frown$ defines a natural equivalence, the
functors $B\to H^k(Y;B)$ also commute with direct limits. Hence by
\cite[Theorem 1]{brown}, $C_*(Y)$ is homotopy equivalent to a complex
$C'$ of f.g.\ (finitely generated) projective modules. Since moreover the
cohomology groups all vanish in dimensions exceeding $n$, we may suppose
by \cite[Theorem E]{FC1} that $C'_r=0$ except when $0\leq r\leq n$.

\begin{proposition}\label{fincond}
Let $Y$ be a connected $\mathit{PD}^n$ complex. Then

\begin{enumerate}
\item[{\rm(i)}] The chain complex $C_*(Y)$ is chain homotopy equivalent to a
complex of f.g.\ projective $R$--modules, vanishing except in dimensions
$r$ with $0\leq r\leq n$.
\item[{\rm(ii)}] The fundamental group $G$ is f.g.\ and a.f.p.\ (almost finitely
presented).
\item[{\rm(iii)}] $Y$ is dominated by a finite complex if and only if $G$
is f.p.
\end{enumerate}
\end{proposition}
We have already proved (i). Since we can attach cells
of dimension $\geq 3$ to $Y$ to obtain a classifying space for $G$,
we can take $C'_2\to C'_1\to C'_0$ as the beginning of a resolution of
$\BZ$ over $R$. Since $C'_1$ is f.g., so is $G$; since $C'_2$ is, $G$
is a.f.p. Now (iii) follows from \cite[Theorem A]{FC1}.

A $\mathit{PD}^n$ complex (or Poincar\'e complex) in general is a complex with
a finite number of components, each of which is a connected $\mathit{PD}^n$
complex. We say that the group $G$ is a $\mathit{PD}^n$ group if $K(G,1)$ is a
$\mathit{PD}^n$ complex: we then have $\cd G=n$.

Corresponding to a manifold with boundary, a connected $\mathit{PD}^n$ pair is a CW pair
$(Y,X)$ with $Y$ connected, with a homomorphism $w\co G\to\{\pm 1\}$ and a class
$[Y]\in {}^tH_n(Y,X;\BZ)$ such that cap product with $[Y]$ induces isomorphisms
$H^r(Y;R) \lraw {}^tH_{n-r}(Y,X;R)$ and $X$ is a $\mathit{PD}^{n-1}$ complex
with fundamental class $\pd_*[Y]$ (so in particular, $w$ induces the
homomorphisms $w$ for the components of $X$). It thus follows from the
five lemma that we have induced isomorphisms $H^r(Y,X;R) \lraw
{}^tH_{n-r}(Y;R)$, and now as before that $[Y]\frown$ induces isomorphisms
$H^i(Y;B)\to {}^tH_{n-i}(Y,X;B)$, $H^i(Y,X;B)\to {}^tH_{n-i}(Y;B)$, for
any $R$-module $B$. If $X=\emptyset$ this reduces to the definition of
$\mathit{PD}^n$ complex. The same arguments as above give:
\begin{addendum}\label{fincondrel}
The conclusions of Proposition~\ref{fincond} apply also if $(Y,X)$
is a $\mathit{PD}^n$ pair.
\end{addendum}
Moreover, since we can find a chain complex for $(\tl{Y},\tl{X})$ with
no 0--cells, it follows from duality that there is a chain complex for
$\tl{Y}$ with no $n$--cells.

Attaching two $\mathit{PD}^n$ pairs by identifying some components of the
boundary yields another such pair. Conversely, there are also results
about cutting along an embedded $\mathit{PD}^{n-1}$ complex. Here we only need
the following \cite[Theorem 2.4]{PC1}.
\begin{proposition}\label{discthm}
Let $(Y,X)$ be a $\mathit{PD}^n$ pair with $n\geq 3$. Then there exists a pair
$(Y',X)$ with $Y'$ dominated by an $(n-1)$--dimensional complex, a map
$f\co S^{n-1}\to Y'$, and a homotopy equivalence $Y'\cup_f e^n\to Y$ (rel
$X$). The triple $(Y',X,f)$ is unique up to homotopy and orientation. If
we suppose (as we may) $f$ an inclusion, $(Y',X\cap S^{n-1})$ is a
$\mathit{PD}^n$ pair.
\end{proposition}
\begin{proof}
In the case when $X$ is empty, this agrees with the result quoted. In
general we first apply the same result to obtain a homotopy equivalence
$h\co X'\cup_g e^{n-1}\to X$. Attach two copies of the mapping cylinder
of $h$ in turn to $Y$ along $X$. The result contains a cell of the form
$e^{n-1}\times e^1$, and it suffices to remove an embedded $n$--cell
from the interior of this to give the existence statement.

If $n=3$ we cannot apply the result to the boundary, but we can use
Theorem~\ref{PD2} below instead.

The same argument as in the case $X$ empty establishes uniqueness in
the general case also.
\end{proof}
\begin{corollary}
The connected sum operation is well defined on connected $\mathit{PD}^n$ pairs
(in the sense that if both are orientable there is a unique connected
sum preserving orientation).
\end{corollary}

A connected $\mathit{PD}^n$ complex or pair with $n\leq 1$ is easily shown to be
homotopy equivalent to a manifold pair: a point, circle or interval. The
same is also known in dimension 2.

\begin{theorem}\label{PD2}
A connected $\mathit{PD}^2$ complex or pair is homotopy equivalent to a compact
manifold pair.
\end{theorem}
The original result was obtained by Eckmann, M\"uller and Linnell in
\cite{eckmul} and \cite{ecklin}. Even stronger results --- in particular,
an analysis of the case when Poincar\'e duality holds over a ring of
coefficients other than $\BZ$ --- are obtained by Bowditch \cite{bowd},
without assuming $G$ f.p.

\section{Decompositions by spheres}\label{secsphere}

To simplify the discussion (and the notation), we restrict from now on to
the orientable case, though the characterisation and splitting results
below were obtained without this restriction. We also assume throughout
that $G$ is f.p., though much is valid without needing this.

Let $(Y,X)$ be a connected $\mathit{PD}^3$ pair. We can define $\what{Y}$
by attaching a 3--disc to $Y$ along each 2--sphere boundary
component. Conversely, we can regard $Y$ as the connected sum of
$\what{Y}$ with a collection of discs $D^3$. We thus suppose from now
on that no component of $X$ is a 2--sphere.

Following 3--manifold terminology, we call a component $X_r$ of $X$
incompressible if the natural map $\pi_1(X_r)\to\pi_1(Y)$ is injective;
$X$ is incompressible if each component is.

Denote the universal cover of $Y$, and the induced coverings of $X$
and its components by adding a tilde.

\begin{lemma}\label{Gfin}
If $Y$ is orientable and $G$ is finite, $X=\emptyset$ and $\tl{Y}\sim S^3$.
\end{lemma}
\begin{proof}
The exact sequence $H_2(Y,X;\BQ)\to H_1(X;\BQ)\to H_1(Y;\BQ)$ is
self-dual, so the image of the first map is a Lagrangian subspace,
of half the dimension. Since $H_1(Y;\BQ)=0$, each component of $X$
has vanishing first Betti number, so is a sphere; hence there are no
components. Now $H_2(\tl{Y};\BZ)=H^1(\tl{Y};\BZ)=0$. The result follows.
\end{proof}

In the nonorientable case one may also have $P^2(\BR)\times I$, for example.

\begin{proposition}\label{ends}
If $Y$ is orientable, the following are equivalent:

\begin{enumerate}
\item[{\rm(i)}] $H^1(G;\BZ G)=0$;
\item[{\rm(ii)}] $e(G)\leq 1$;
\item[{\rm(iii)}] $e(\tl{Y})\leq 1$;
\item[{\rm(iv)}] $\pi_2(Y)=0$ and $X$ is incompressible.
\end{enumerate}
\end{proposition}
\begin{proof}
Here $\BZ$ can be replaced by a field $k$.

Since we can obtain a $K(G,1)$ by attaching cells of dimension $\geq 3$
to $Y$, $H^1(G;\BZ G)=H^1(Y;\BZ G)=H^1_c(\tl{Y})$ is the cohomology
of $\tl{Y}$ calculated with finite cochains. Denote by $C^*(\tl{Y})$,
$C^*_c(\tl{Y})$ the chain complex of cochains of $\tl{Y}$ and the
subcomplex of finite cochains; write $C^*_e(\tl{Y})$ for the quotient,
and $H^*_e(\tl{Y})$ for its cohomology groups. From the exact sequence
\begin{equation}\label{endseq}
 H^0_c(\tl{Y})\to H^0(\tl{Y})\to H^0_e(\tl{Y})\to H^1_c(\tl{Y})\to
 H^1(\tl{Y}),
\end{equation}
where $e(G)=e(\tl{Y})$ is the rank of the middle term, we see first that
if $G$ is finite $e(G)=0$ and $H^1_c(\tl{Y})=0$; thus by the lemma, all
of (i)--(iv) hold. From now on, assume $G$ infinite. Then the extreme
terms of (\ref{endseq}) vanish and $H^0(\tl{Y})=\BZ$, so $e(G)\geq
1$, with equality holding only if $H^1_c(\tl{Y})=0$.  Thus (i)--(iii)
are equivalent.

Observe that, for each component $X_r$ of $X$, the image of
$\pi_1(X_r)\to\pi_1(Y)$ is infinite. For, if not, the map $H_1(X_r;\BQ)\to
H_1(Y;\BQ)$ would be zero. This contradicts the facts that the kernel of
$H_1(X;\BQ)\to H_1(Y;\BQ)$ is a Lagrangian subspace and, since $X_r$ is
not a sphere, $H_1(X_r;\BQ)$ has non-vanishing intersection numbers. Hence
the kernel, $J_r$ say, of $\pi_1(X_r)\to\pi_1(Y)$ has infinite index
in $\pi_1(X_r)$, so is a free group. Then each component of $\tl{X_r}$
is non-compact, and $H_1(\tl{X_r};\BZ)$ is a sum of copies of $J_r^\mathit{ab}$.

Now by duality $H^1_c(\tl{Y};\BZ)\cong H_2(\tl{Y},\tl{X};\BZ)$,
and in the exact sequence \[H_2(\tl{X};\BZ)\to H_2(\tl{Y};\BZ)\to
H_2(\tl{Y},\tl{X};\BZ)\to H_1(\tl{X};\BZ)\to H_1(\tl{Y};\BZ),\] the
extreme terms vanish: the left hand one since each component of $\tl{X}$
is non-compact, the right since we have the universal cover. Hence
$H_2(\tl{Y},\tl{X};\BZ)$ vanishes if and only if (a) $\pi_2(Y)\cong
H_2(\tl{Y};\BZ)$ vanishes and (b) $H_1(\tl{X};\BZ)$ vanishes, ie each
$J_r^\mathit{ab}$ does. But since $J_r$ is free, this implies that $J_r$ is
trivial, hence $X_r$ incompressible.
\end{proof}

In the case when $Y$ is a closed orientable 3--manifold, the sphere
theorem states that the following are equivalent:

\begin{enumerate}
\item[(i)] $e(G)\geq 2$ and $G\not\cong\BZ$,
\item[(ii)] $G$ is a free product,
\item[(iii)] $Y$ splits non-trivially as a connected sum.
\end{enumerate}
Moreover, if $G\cong \BZ$, $Y$ is homeomorphic to $S^2\times S^1$.

If $G\not\cong\BZ$ and $e(G)\geq 2$, $G$ is a free product, and we can
decompose $Y$ as a connected sum. By Gru\v{s}ko's Theorem, the process
of consecutive decompositions must terminate. We end with $G$ being
a free product of free factors each of which is either ($e=2$) $\BZ$,
($e=0$) finite, or ($e=1$) the fundamental group of a 3--manifold which,
by Proposition~\ref{ends}, is aspherical, so has torsion-free fundamental
group. Thus $G$ is a free product of finite groups and a torsion-free
group. It follows that any finite subgroup of $G$ is contained in a
free factor.

It is natural to hope for corresponding results for orientable
$\mathit{PD}^3$
complexes: we shall see that this is too optimistic. We still know,
by \cite{PC1}, that if $G\cong \BZ$, $Y$ is homotopy equivalent to
$S^2\times S^1$. Next there is a recent result of Crisp.

\begin{theorem}\label{crisp} {\rm\cite{crisp}}\qua
If $Y$ is an orientable $\mathit{PD}^3$ complex and $e(G)\geq 2$, then either
(a) $G$ is a free product or (b) $G$ is virtually free.
\end{theorem}
Since $e(G)\geq 2$, it follows from Stallings' theorem that there is a non-trivial
action of $G$ on a tree $T$ with finite edge groups $G_e$; moreover, as $G$ is
accessible, we may suppose that each vertex group $G_v$ has at most 1 end. In the
exact sequence
\begin{multline*}
0\to H^0(G;M)\to\oplus_vH^0(G_v;M)\to \\
 \oplus_eH^0(G_e;M)\to H^1(G;M)\to\oplus_v H^1(G_v;M)
\end{multline*}
take $M=\BZ G$. Since, for any subgroup $K\subset G$, $H^i(K;\BZ
G)=H^i(K;\BZ K)\otimes_{\BZ K}\BZ G$, and for $K$ infinite and 1--ended,
$H^0(K;\BZ K)=H^1(K;\BZ K)=0$, while for $K$ finite $H^0(K;\BZ K)=\BZ$
and $H^1(K;\BZ K)=0$, the sequence reduces to
$$0\to \oplus_{G_v\,\text{finite}}(\BZ\otimes_{\BZ G_v}\BZ G)\to
\oplus_e(\BZ\otimes_{\BZ G_e}\BZ G)\to H^1(G;\BZ G)\to 0.$$
Now Crisp regards $H^1(G;\BZ G)$ as a modification $\Pi(T)$ of $H^1_c(T)$
in which the vertices with infinite stabilisers are omitted from the
calculation, and proves that $\Pi(T)$ is free over $\BZ$ with rank
$\max(0,e(T)+\infty(T)-1)$, where $e(T)$ is the number of ends of $T$
and $\infty(T)$ is the number of vertices with infinite stabilisers.

If any edge group is trivial, $G$ splits as a free product, and if all
vertex groups are finite, $G$ is virtually free. The key idea of the
proof is to consider a finite cyclic subgroup $C$ of an edge group and
compare a calculation $H_s(C;H^1(G;\BZ G))\cong H_{s+3}(C;\BZ)$ using
Poincar\'e duality with the above.

It is natural to expect that case (b) cannot occur unless $G$ itself is free. An
analysis by Hillman showed that the first case which cannot be resolved by easy
arguments is the free product of two copies of the dihedral group of order 6,
amalgamated along a subgroup of order 2. In a recent preprint he established the
following, thus showing that case (b) does indeed occur.
\begin{theorem}\label{counterex}{\rm\cite{hillmanpre}}\qua
There is an orientable $\mathit{PD}^3$ complex whose fundamental group is the
amalgamated free product $D_6 *_{\BZ_2} D_6$.
\end{theorem}
The proof depends on the criterion of Turaev to be discussed below
(Theorem~\ref{turaevcharsn}).

The first result indicating a purely algebraic treatment of $\mathit{PD}^3$
complexes was the following theorem of Hendriks:
\begin{proposition}\label{classify} {\rm\cite{hendriks}}\qua
Two Poincar\'e 3--complexes with the same fundamental group $G$ are
homotopy equivalent if and only if the images of their fundamental
classes in $H_3(G;\BZ)$ coincide.
\end{proposition}
Next Turaev in \cite{tur1} obtained a characterisation of which pairs
$(G,z)$ with $z\in H_3(G;\BZ)$ correspond to $\mathit{PD}^3$ complexes. In
\cite{tur2} he used this, together with Hendriks' theorem, so show that
if the fundamental group of a $\mathit{PD}^3$ complex splits as a free product,
there is a corresponding split of the complex as a connected sum. In a
subsequent paper \cite{tur3} he gave an improved and unified version of
all three results.

We now review Turaev's argument in our own notation. Let $\mathit{Mod}_R$
denote the category of (f.g.\ right) $R$--modules; define a morphism to
be nullhomotopic if it factors through a projective module. Form the
quotient category $\mathit{PMod}_R$ by these morphisms, and write $[M,N]$
for the group of morphisms from $M$ to $N$ in $\mathit{PMod}_R$. Equivalence
in $\mathit{PMod}_R$ is known as {\em stable equivalence}; $R$--modules $M,M'$
are stably equivalent if and only if there exist f.g.\ projective modules
$P,P'$ such that $M\oplus P$ and $M'\oplus P'$ are isomorphic.

We introduce abbreviated notation as follows. Write $K:=K(G,1)$: consider
spaces mapped to $K$ and for each such space $W$, write $\wtl{W}$ for
the covering space of $W$ induced from the universal cover of $K$, and
$C_*(W)$, $C^*(W)$ for the chain and finite cochain groups of $\wtl{W}$,
considered as free modules over $R:=\BZ G$. We use a corresponding
notation for pairs. Also where coefficients for homology or cohomology
are unspecified, they should be understood as $R$. We write $\fG:=\Ker
(R\to\BZ)$ for the augmentation ideal of $G$.

The following algebraic construction is the key to the argument. For
any projective chain complex $C_*$ over $R=\BZ G$, write $IC_*$
for $\fG\otimes_R C_*$ and set $F^rC:=\Coker(\gd^{r-1}\co C^{r-1}\to
C^r)$. Then there is a homomorphism $\nu_{C,r}\co H_r(IC_*)\to [F^rC,\fG]$
induced by evaluating a representative cocycle in $C^r$ (for an element of
$F^rC$) on a representative cycle in $\fG\otimes_R C_r$ (for an element
in $H_r(IC_*)$). This construction is natural with respect to maps of
chain complexes $C_*$. Moreover, Turaev shows that:
\begin{lemma}\label{turaevlem}
For any projective chain complex $C_*$, $\nu_{C,r}\co H_r(IC_*)\to
[F^rC,\fG]$ is an isomorphism.
\end{lemma}

If $C_*$ is a positive complex with $H_0(C)\cong\BZ$, then the image of
$d_1\co C_1\to C_0$ is stably equivalent to $\fG\lhd\BZ G$ so that if also
$H_1(C)=0$, $\Coker d_2\co C_2\to C_1$ is stably equivalent to $\fG$. A
class $\mu\in H_{r+1}(C;\BZ)$ induces a chain mapping $C^{r+1-*}(K)\to
C_*(K)$, which is determined up to chain homotopy. There is an induced
map of the cokernel $F^rC$ of $\gd\co C^{r-1}(K)\to C^r(K)$ to the
cokernel of $d\co C_2(K)\to C_1(K)$, and hence to the augmentation ideal
$\fG$. Denote the composite by $\nu(\mu)\in [F^r(C),\fG]$.

We also have an exact sequence
$$H_{r+1}(C)\to H_{r+1}(C;\BZ)\stackrel{\pd_*}{\lraw} H_r(IC_*)\to H_r(C).$$
Turaev also shows that $\nu(\mu)=\nu_{C,r}(\pd_*\mu)$. In particular,
if also $H_{r+1}(C)=H_r(C)=0$, then $\pd_*$ is an isomorphism, so $\nu_{C,r}$
induces an isomorphism $\nu\co H_{r+1}(C;\BZ)\to [F^rC,\fG]$.

For our purposes we need only the case $r=2$, so omit $r$ from the notation.

If $Y$ is an oriented $\mathit{PD}^3$ complex with fundamental class $[Y]\in
H_3(Y;\BZ)$, we have a classifying map $i\co Y\to K$ and hence a class
$i_*(Y)\in H_3(G;\BZ)$. Now Turaev's main theorem (in the orientable
case) is:
\begin{theorem}\label{turaevcharsn}
Given a triple $(G,\mu)$ with $G$ an (f.p.) group and $\mu\in H_3(G;\BZ)$,
there is an oriented $\mathit{PD}^3$ complex $Y$ with fundamental group $G$
and $i_*[Y]=\mu$ if and only if $\nu(\mu)$ is a stable equivalence.

Moreover, if this holds, $Y$ is unique up to oriented homotopy equivalence.
\end{theorem}
Necessity of the condition follows easily from the definition of
$\nu$. The key step in the proof of sufficiency is the construction of
the complex $Y$.

Let $Z$ be a finite connected 2--complex with $\pi_1(Z)=G$, eg the
2--skeleton of $K$.  Our hypothesis gives a preferred stable equivalence
from the cokernel $F^2C(Z)=C^2(Z)/\gd(C^1(Z)$ to $\fG$. Replacing $Z$,
if necessary, by its bouquet with a number of 2--spheres, we may suppose
given an isomorphism $\phi$ of the cokernel $F^2C(Z)$ to $\fG\oplus P$
for some f.g.\ projective $R$--module $P$.

We can make $P$ free as follows. By a lemma of Kaplansky, for any
f.g.\ projective module $P'$ there exist (infinitely generated) free
modules $F,F'$ with $F'\cong F\oplus P'$. We can thus attach 2--spheres
to $Z$ corresponding to generators of $F'$ and 3--cells to kill the
generators of $F$: this has the same effect as adding a copy of $P'$ to
$C_2(Z)$, while the new complex $Z$ still has the homological properties
of a 2--dimensional complex. If we choose $P'$ appropriately, this will
make $P$ free and f.g., of rank $t$, say.

Write $d^3$ for the map $C^2(Z)\to R^{t+1}$ defined by composing $\phi$
with the projection. Dualising gives a homomorphism $d_3\co R^{t+1}\to
C_2(Z)$ such that $d_2\circ d_3=0$. Thus the image of $d_3$ is contained
in $Z_2(Z)= \pi_2(\wtl{Z})=\pi_2(Z)$.  Attach 3--cells to $Z$ by maps
$S^2\to Z$ corresponding to the images of the generators of $R^{t+1}$:
this gives a complex $Y$.

We now show that this complex $Y$ is a Poincar\'e 3--complex:
first we calculate $H_3(Y;\BZ)$. The map $d^3$ is the composite
$C^2(Z)\twoheadrightarrow\fG\oplus R^t\hookrightarrow R^{t+1}$,
which is the direct sum of an isomorphism on $R^t$ and a composite
$A\twoheadrightarrow\fG\hookrightarrow R$. When we tensor over $R$
with $\BZ$, the map $\fG\hookrightarrow R$ gives 0. Dually, $d_3$
with coefficients $\BZ$ is the direct sum of an isomorphism of $\BZ^t$
on itself and a zero map of $\BZ$.  Thus $H_3(Y;\BZ)\cong\BZ$. Also,
there is a preferred generator, giving a class $[Y]$.

We next show that $i_*[Y]$ is equal to the given class $\mu$. For this
we apply Turaev's Lemma~\ref{turaevlem} to $C_*(K)$. Since $\wtl{K}$
is contractible, the hypothesis $H_3(K,X)=H_2(K,X)=0$ is satisfied. By
construction, $\nu[Y]$, which by naturality is equal to $\nu i_*[Y]$,
is the stable isomorphism $\nu(\mu)$. It follows that indeed $i_*[Y]=\mu$.

The map $\phi_*$ is part of an exact triangle $C^{3-*}(Y)\to C_*(Y)\to
D_*\to C^{4-*}(Y)$ in the derived category. Since $\phi_*$ is self-dual,
so is $D_*$. By construction, $H^2(Y)=0$ and $H^3(Y)\cong\BZ$, and
so $\phi_*$ induces isomorphisms $H^2(Y)\to H_1(Y)$ and $H^3(Y)\to
H_0(Y)$. In particular, $H_0(D)\cong H_1(D)=0$, so $D$ is chain equivalent
to a complex $0\to D_4\to D_3\to D_2\to 0$. As $D_*$ is self-dual, we may
similarly eliminate $D_4$ and $D_3$, leaving only a projective module
$D_2'$, say. Thus $C_2\to D'_2$ is a split surjection and $H_2(Y)\cong
D'_2\oplus H^1(Y)$. By duality, $D'_2\to C^2$ is a split injection
and $H^2(Y)\cong D'_2\oplus H_1(Y)$. As $H^2(Y)=0$, $D'_2=0$ and $D_*$
is acyclic.

As to uniqueness, we  may construct a $K(G,1)$ complex $K$ by attaching
cells of dimension $\geq 3$ to $Y$. If $Y'$ is another $\mathit{PD}^3$ complex with
the same fundamental group $G$, there is no obstruction to deforming
a map $Y'\to K$ inducing an isomorphism of $G$ to a map into $Y$;
moreover if we assume that the fundamental classes have the same image
in $H_3(K;\BZ)$, a careful argument shows that we may suppose $Y'\to Y$
of degree 1. It is now easy to see that we have a homotopy equivalence.

\begin{theorem}\label{turaevsplit} {\rm\cite{tur2}}\qua
If $Y$ is a $\mathit{PD}^3$ complex such that $G=\pi_1(Y)=G'*G''$ is a free
product, $Y$ is homotopy equivalent to the connected sum of $\mathit{PD}^3$
complexes $Y'$ and $Y''$, with $\pi_1(Y')\cong G'$ and $\pi_1(Y'')\cong
G''$.
\end{theorem}
The image of the fundamental class of $Y$ gives an element $\mu\in
H_3(G;\BZ)\cong H_3(G';\BZ)\oplus H_3(G'';\BZ)$, and hence classes
$\mu'\in H_3(G';\BZ)$, $\mu''\in H_3(G'';\BZ)$. It will suffice
by Theorem~\ref{turaevcharsn} to show that $\nu(\mu')$ is a stable
equivalence.

Choose Eilenberg--MacLane 2--complexes $K',K''$, with respective
fundamental groups $G',G''$ and group rings $R',R''$; then $K=K'\vee K''$
is a $K(G,1)$. We have $C_2(K')\to C_1(K')\to \fG'\to 0$, and similarly
for $G''$. Tensor these over $R',R''$ with $R$ and add. We obtain
$$C_2(K)\to C_1(K) \to (\fG'\otimes_{R'}R)\oplus (\fG''\otimes_{R''}R) \to 0.$$
Thus $(\fG'\otimes_{R'}R)\oplus (\fG''\otimes_{R''}R)\cong \fG$.

We have $\nu(\mu)\in [F^2C(K),\fG]$, and similarly for $G'$, $G''$. Since
$C^*(K)\cong (C^*(K')\otimes_{R'}R)\oplus (C^*(K'')\otimes_{R''}R)$,
we can identify $F^2C(K)$ with $(F^2C(K')\otimes_{R'}R)\oplus
(F^2C(K'')\otimes_{R''}R)$. Indeed, with the obvious interpretation,
we can write $\nu(\mu)=(\nu(\mu')\otimes_{R'}R)\oplus
(\nu(\mu'')\otimes_{R''}R)$.

We now apply $\otimes_R R'$. If $M$ denotes either $\fG$ or $F^2C(K)$,
we can write $M\cong (M'\otimes_{R'}R)\oplus (M''\otimes_{R''}R)$. The
first summand yields $(M'\otimes_{R'}R)\otimes_R R'=M'$; the second gives
$$(M''\otimes_{R''}R)\otimes_R R'=(M''\otimes_{R''}\BZ)\otimes_{\BZ} R'.$$
In each case, $M''$ is f.g.\ over $R''$, so $M''\otimes_{R''}\BZ$ is
an f.g.\ abelian group, which we can express as $F\oplus T$ with $F$
free abelian and $T$ finite abelian. The summands $F\otimes_{\BZ}R'$
can be ignored, as we are only concerned with stable isomorphism.

It follows that the given stable isomorphism $\nu(\mu)$ gives,
on tensoring over $R$ with $R'$, the direct sum of $\nu(\mu')$,
which is a map of torsion-free modules, and a map of torsion modules
$T_1\otimes_{\BZ}R' \to T_2\otimes_{\BZ}R'$, which must thus be an
isomorphism. Hence indeed $\nu(\mu')$ is a stable isomorphism.

\section{Decomposition by spheres and discs}\label{secdisc}

In the case $X\neq\emptyset$, matters are distinctly more
complicated. Lemma 2.1, Proposition 2.2 and Theorem 2.3 were already
framed to include the general case. We now seek to decompose $(Y,X)$ in
some way (some analogue of a decomposition of a 3--manifold by embedded
spheres and discs) until the conditions of Proposition 2.2 hold for each
piece $(Y,X)$ of the decomposition. Then for a piece with fundamental
group $G$, $G$ is not a free product and either
\begin{enumerate}
\item[(i)] $e(G)=0$, $G$ is finite, $Y$ is finitely covered by a homotopy $S^3$,
\item[(ii)] $e(G)=2$ and $(Y,X)$ is one of few possibilities,
\item[(iii)] $e(G)=1$, $X$ is incompressible and $Y$ aspherical.
\end{enumerate}

Although this is overoptimistic, we will investigate how far one
can go towards such a result: our overall conclusion is that results
corresponding to the manifold case can be proved if either $G$ is torsion
free or $X$ is incompressible.

If $e(G)<2$, we can apply Proposition~\ref{ends}, so $\pi_2(Y)=0$. If
$e(G)=1$, $Y$ is aspherical: we treat this case in the next section. If
$e(G)=0$, $G$ is finite, so $X=\emptyset$, and $\tl{Y}$ is homotopy
equivalent to $S^3$. For the case $e(G)=2$, we have:
\begin{proposition}\label{2ends}
An orientable $\mathit{PD}^3$ pair $(Y,X)$ such that $e(G)=2$ is homotopy
equivalent to one of $(P^3(\BR)\# P^3(\BR),\emptyset)$, $(S^2\times
S^1,\emptyset)$ and $(D^2\times S^1,S^1\times S^1)$.
\end{proposition}
\begin{proof}
If $X=\emptyset$, then it follows from \cite[Theorem 4.4]{PC1} that either
$Y\simeq P^3(\BR)\# P^3(\BR)$ (so $G$ is a free product) or $Y\simeq
S^2\times S^1$. In general $G$ has a finite normal subgroup $F$ with
quotient either $\BZ_2 *\BZ_2$ or $\BZ$, so $H_1(Y;\BZ)\cong G^\mathit{ab}$
has rank 0 or 1. Thus $H_1(X;\BZ)$ has rank at most 0 resp. 2. We have
already dealt with the case $X=\emptyset$, so may suppose $G/F\cong\BZ$
and $X$ a torus.

It follows from (\ref{endseq}) that $H^1_c(\tl{Y};\BZ)\cong\BZ$,
so by duality $H_2(\tl{Y},\tl{X};\BZ)\cong\BZ$. Since the image of
$H_1(X;\BZ)\to H_1(Y;\BZ)$ has rank 1, each component of $\tl{X}$ is
homotopy equivalent to $S^1$. It now follows from the exact sequence
$$H_2(\tl{X};\BZ)\to H_2(\tl{Y};\BZ)\to H_2(\tl{Y},\tl{X};\BZ)\to
  H_1(\tl{X};\BZ) \to H_1(\tl{Y};\BZ),$$
whose end terms vanish, that $H_2(\tl{Y};\BZ)=0$, so $Y$ is
aspherical. Hence $G$ is torsion free, so $G\cong\BZ$. Moreover, $\tl{X}$ is
connected, so $\pi_1(X)\to\pi_1(Y)$ is surjective. Thus indeed
$(Y,X)\simeq (D^2\times S^1,S^1\times S^1)$.
\end{proof}

Since $e(G)\geq 2$ did not imply $G$ a free product in the preceding
section, we cannot hope to do better here. In some sense, things are now
no worse. For suppose $(Y,X)$ an orientable $\mathit{PD}^3$ pair. Form the double
$DY$, with fundamental group $\what{G}$, say. Suppose $\what{G}$ is a
free product of finite groups and a torsion free group. By Kuro\v{s}'
subgroup theorem, the same follows for any subgroup of $\what{G}$. But
the natural map $G\to\what{G}$ is injective since there is a retraction
by folding the factors of the double. Thus any finite subgroup of $G$
is contained in a free factor, and Theorem 2.3 now shows that $G$ is a
free product of the desired type.

Crisp, in \cite{crisp}, gives an extension of Theorem~\ref{crisp} to
the case of Poincar\'e pairs $(Y,X)$. First he shows that his argument
remains valid if $X$ is incompressible.
\begin{theorem}\label{crisp2}
If $(Y,X)$ is an orientable $\mathit{PD}^3$ pair with $X$ incompressible and
$e(G)\geq 2$, then either (a) $G$ is a free product or (b) $G$ is
virtually free and $X=\emptyset$.
\end{theorem}
It remains only to observe in the second case that a non-trivial
fundamental group of a closed surface is never virtually free.

Next Crisp observes that if the loop theorem is applicable, he can reduce
to the incompressible case by attaching 2--handles to $Y$. We would thus
like a version of the loop theorem, and next digress to discuss this.

The original Loop Theorem for 3--manifolds was proved by Papakyriakopoulos
\cite{papa}, further proofs were given by Stallings \cite{stall} and
Maskit \cite{maskit}. The result we would like (shorn of refinements) is:
\begin{conjecture}\label{PDloopthm}
Let $(Y,X)$ be a $\mathit{PD}^3$ pair with $X$ a 2--manifold; let $F$ be a
component of $X$ such that $\pi_1(F)\to\pi_1(Y)$ is not injective. Then
there exists a simple loop $u$ in $F$ which is nullhomotopic in $Y$
but not in $F$.
%and $N\lhd\pi_1(F)$ a normal subgroup such that
%$\Ker(\pi_1(F)\to\pi_1(Y))\not\subseteq N$. Then there exists a simple loop $u$ in
%$F$, nullhomotopic in $Y$, whose homotopy class is not in $N$.
\end{conjecture}
This result is claimed in \cite{thomas}, but the proposed proof appears
to have gaps.  The only other relevant reference known to the author is
Casson--Gordon \cite{cassgor}, but (despite the assertion in the review
\href{http://www.ams.org/mathscinet-getitem?mr=0722728}{MR0722728} 
that it ``proves, via a geometric argument on planar coverings,
that the loop theorem of Papakyriakopoulos is true for surfaces that
bound mere 3--dimensional duality spaces") the result obtained in that
paper is of a somewhat different nature, and works only modulo the
intersection of the dimension subgroups of $G$.

We obtain a partial result by following Maskit's proof of the theorem
for 3--manifolds.
\begin{lemma}
Let $(Y,X)$ be a $\mathit{PD}^3$ pair with $X$ a 2--manifold, with universal
covering $\tl{Y}$ and induced covering $\tl{X}$. Then every component
of $\tl{X}$ is planar (ie can be embedded in a plane).
\end{lemma}\begin{proof}
Suppose, if possible, there is a non-planar component $\tl{F}$. Since
$\tl{Y}$ and hence $\tl{X}$ is orientable, it follows from the standard
theory of surfaces that there exist simple loops $A,B\subset\tl{F}$ which
meet transversely in just one point. As $\tl{Y}$ is simply-connected,
$B$ bounds a 2--cycle, which thus defines a homology class in
$H_2(\tl{Y},\tl{X};\BZ)$. This has a dual cohomology class $\gb\in
H^1_c(\tl{Y};\BZ)$.

Write $i\co\tl{X}\to\tl{Y}$ for the inclusion, and $\ga$ for the class
of $A$ in $H_1(\tl{X};\BZ)$. It follows from the definition of $\gb$
that $i^*\gb$ is the cohomology class dual to the cycle $B$. From
the choice of $A$ and $B$, $\langle i^*\gb,\ga\rangle=1$. Hence
\[\langle\gb,i_*\ga\rangle=\langle i^*\gb,\ga\rangle=1.\] Thus $0\neq
i_*\ga\in H_1(\tl{Y};\BZ)$, contradicting simple connectivity of $\tl{Y}$.
\end{proof}
\begin{proposition}\label{planarity} {\rm\cite[Theorem 3]{maskit}}\qua
Let $p\co\tl{F}\to F$ be a regular covering with $F$ a compact surface
and $\tl{F}$ planar. Then there exist a finite disjoint set of simple,
orientation-preserving loops $u_1,\ldots,u_s$ on $F$ and positive integers
$n_1,\ldots,n_s$ such that $\pi_1(\tl{F})$ is the normal subgroup of
$\pi_1(F)$ generated by $u_1^{n_1},\ldots,u_s^{n_s}$.
\end{proposition}
Applying this to the covering of $F$ induced by the universal covering
of $Y$ in Hope~\ref{PDloopthm}, which is planar by the lemma, gives
simple loops $u_i$ in $F$ with $u_i^{n_i}$ nullhomotopic in $Y$. In
the manifold case it can now be shown (see \cite[p 287]{papa} that for
each such simple loop $u_i$, supposed orientation-preserving, $u_i$ is
already nullhomotopic in $Y$. Here, at least if $\pi_1(Y)$ is torsion
free, it follows that each $u_i$ is already nullhomotopic in $Y$. Thus:
\begin{theorem}\label{torfreeloop}
Hope~\ref{PDloopthm} is true if $G$ is torsion free.
\end{theorem}
This is of no use for applying Crisp's theorem since if $e(G)\geq 2$
and $G$ is torsion free, it follows anyway that either $G\cong\BZ$ or $G$
is a free product.

Now suppose that $G$ splits as a free product: can one establish
the existence of some kind of splitting of $Y$? Let us begin by
considering what happens in the case of manifolds. Suppose given a
splitting $G=G'*G''$. Join the base points of $K'$ and $K''$ by an
arc to form the Stallings wedge $K$ (which is indeed a $K(G,1)$), take
the induced map $h\co Y\to K$, and make $h$ transverse to the centre
point of the arc.  The pre-image is an orientable surface $V$ with each
boundary loop null-homotopic in $Y$. For each component $V_i$ of $V$,
$\pi_1(V_i)$ maps to 0 in $G$, so by the loop theorem, if $V_i$ is not
simply connected, there is a compressing disc. An argument of Stallings
now shows that any surgery on $V$ performed using a compressing disc
can be induced by a homotopy of $h$. We thus reduce to the case when
each $V_i$ is simply-connected, hence a sphere or disc.

Consider a single $V_i=V$, cut $Y$ along it to give $Y^0$, and attach
3--discs to the new boundary components to obtain $\what{Y^0}$. First
suppose $V$ is a sphere. If $V$ separates $Y$, $Y$ is a connected sum
of the two components of $\what{Y^0}$. If one of these has trivial
fundamental group, it is a homotopy sphere and we can deform $h$ to
remove the component $V$. If $V$ fails to separate $Y$, the union $Y'$
of a collar neighbourhood of $V$ and the neighbourhood in $Y$ of an
embedded circle meeting $V$ transversely in one point has boundary a
disc and closed complement $Y''$, say, so $Y$ is a connected sum of
$\what{Y''}$ and $\what{Y'}\cong S^2\times S^1$.

If $V$ is a disc, and separates $Y$, then $Y$ is a boundary-connected
sum of the components of $Y^0$; ie it is formed by identifying 2--discs
embedded in the boundaries of these components. If either of these
components is simply-connected, it is contractible, and we can deform
to remove the component $V$.

If $V$ is a disc which fails to separate $Y$, and $\pd V$ fails to
separate the relevant component $X_r$ of $X$, let $Y'$ be a regular
neighbourhood of the union of $V$ and an embedded circle in $X_r$ meeting
$V$ transversely in one point. Then $Y'$ is homeomorphic to $S^1\times
D^2$ and has relative boundary a 2--disc, so $Y$ is a boundary-connected
sum of $S^1\times D^2$ with some $Y''$; if $Y''$ is simply-connected,
$Y$ itself is homotopy equivalent to $S^1\times D^2$. Finally if $\pd V$
separates $X_r$ and $V$ doesn't separate $Y$, then $Y$ is formed from
$Y^0$ by identifying discs embedded in distinct boundary components ---
a sort of boundary-connected sum of $Y^0$ with itself.

We may thus speculate that if $(Y,X)$ is a $\mathit{PD}^3$ pair and
$G=\pi_1(Y)$ is a non-trivial free product, there is a splitting of
$(Y,X)$ as either connected sum or boundary-connected sum corresponding
to a non-trivial splitting of $G$, or self-boundary-connected sum. This
seems somewhat complicated, so we look for an alternative approach.

In the case when $Y$ is a manifold, the result of cutting $Y$ along
a properly embedded disc $(D^2,S^1)\subset (Y,X)$ is homeomorphic to
that obtained by removing from $Y$ the (relative) interior of a collar
neighbourhood $(D^2,S^1)\times D^1$.  This is homeomorphic to the
result of the following sequence of operations:
\begin{enumerate}
\item[(i)] attach a 2--handle $(D^2,S^1)\times D^1$ to $Y$ using the same
embedding $S^1\times D^1\to X$;
\item[(ii)] cut along the copy of $S^2$ which is the union of the copy
of $D^2\times 0$ embedded in $Y$ and the copy we have just attached;
\item[(iii)] attach a 3--disc to each of the copies of $S^2$ on the boundary
of the result.
\end{enumerate}
Thus instead of cutting along discs we will consider a
sequence of operations of attaching 2--handles, and then splitting along
embedded spheres.

If this sequence of operations leads to a manifold with incompressible
boundary, then this property must already hold at the stage when we have
done the handle attachments. Thus our new plan is: first attach 2--handles
to the $\mathit{PD}^3$ pair $(Y,X)$ to make the boundary incompressible;
then split by embedded 2--spheres as long as the fundamental group is a
free product. To attach 2--handles, we require simple loops in $X$. Since
we wish to apply the loop theorem, let us assume $\pi_1(Y)$ torsion free.
Then if $X_r$ is compressible, there is an essential embedding of $S^1$ in
$X_r$ which is null-homotopic in $Y$. Let $X_b$ be a regular neighbourhood
of the image of $S^1$ (we think of this part of $X$ as coloured black).

Suppose inductively that we have a compact 2--dimensional submanifold
$X_b$ of $X$ such that the composite $X_b\subset X\subset Y$ is
nullhomotopic. Write $X_w$ for the closure of $X\sm X_b$ (the `white
part'). If, for some component $Z$ of $X_w$, $\pi_1(Z)\to \pi_1(Y)$ is
not injective, the loop theorem provides an essential embedding of $S^1$
in $Z$ which is nullhomotopic in $Y$, and we add a neighbourhoood of this
loop to $X_b$. If $Z$ has a boundary, we also add to $X_b$ an arc joining
this neighbourhood to the boundary. We can also attach a 2--handle to $Y$
by the chosen curve at each stage of the construction. These attachments
do not affect $\pi_1(Y)$ since the curves were nullhomotopic. The
boundary of the result is obtained from $X_w$ by attaching a 2--disc
to each boundary component. No component of this can be a 2--sphere,
since this would mean we had used an inessential curve.

At each step we are doing a surgery on the boundary. Thus the set of
genera of the components is changed either by decreasing one by 1 or by
replacing $g_1+g_2$ by $g_1$ and $g_2$. As all genera are positive, the
procedure must terminate. Thus there is a disjoint union of embedded
copies of $S^1$ in $X$ such that, for each component $X_r$ of $X$,
the classes of the circles embedded in $X_r$ generate the kernel of
$\pi_1(X_r)\to\pi_1(Y)$ as a normal subgroup.

At each stage each $X_b\cap X_r$ is connected. Recall that the kernel
of $H_1(X;\BQ)\to H_1(Y;\BQ)$ is a Lagrangian subspace with respect
to the intersection product. Thus for each component $X_r$ of $X$,
$Ker(H_1(X_r;\BQ)\to H_1(Y;\BQ))$ is isotropic. Hence the image of
$H_1(X_b\cap X_r;\BQ)$ in $H_1(X_r;\BQ)$ is isotropic.  It follows that
the surface $X_b\cap X_r$ is planar. Suppose $X_b\cap X_r$ has $m+1$
boundary components. Then our procedure involved $m$ steps acting on
$X_r$. These can be performed along any $m$ of the $m+1$ boundary curves
of $X_b\cap X_r$.

If the result of this construction were unique, it would follow that
we had lost no information by doing all the handle additions first,
since there was only one way to make the boundary incompressible. But
it is not clear that any simple loop in $X_r$ which is nullhomotopic in
$Y$ can be isotoped to lie in $X_b\cap X_r$. Or we can consider taking
$X$, collapsing $X_b$ to a point, and take the map $\pi_1(X/X_b)\to
\pi_1(Y)$. Now $X/X_b$ is a bouquet of closed surfaces, and the
fundamental group of each injects in $G$, but it does not follow that
$\pi_1(X/X_b)$ (which is their free product) does. When this is the case,
the result of the construction is essentially unique.

From now on we consider a (connected) orientable $\mathit{PD}^3$
pair $(Y,X)$ with incompressible boundary. We would like to show
that any splitting of $G$ is induced by a connected sum splitting
of $Y$. Following the arguments of Turaev in the preceding section,
we first seek a characterisation of the possible homotopy types of
$(Y,X)$. Here we follow the preprint \cite{bleile} of Bleile, but we
omit most of the details.

In view of Theorem~\ref{PD2} we may suppose given an oriented 2--manifold
$X$, with no 2--sphere components, a group $G$, and a homotopy class of
maps $X\to K$. We will later assume ($X$ is incompressible) that each map
$H_r\to G$ is injective, thus in this case each component of $\tl{X}$
is contractible. Suppose also given a class $\mu\in H_3(K,X;\BZ)$ such
that $\pd_*(\mu)=[X]$ is the fundamental class of $X$.

We may suppose the map $X\to K$ an inclusion, so have an exact sequence
$C_*(X)\to C_*(K)\to C_*(K,X)$ of $R$--free chain complexes. As before,
a class $\mu\in H_3(K,X;\BZ)$ induces cap products $\mu\frown\co
H^{3-i}(K,X;-)\to H_{i}(K;-)$ via a chain mapping $C^{3-*}(K,X)\to
C_{*}(K)$, which is determined up to chain homotopy.  The induced
map of the cokernel $FC^2(K,X)$ of $\gd\co C^1(K,X)\to C^2(K,X)$
to the cokernel of $d\co C_2(K)\to C_1(K)$ is unique up to stable
equivalence. As before, since $C_*(K)$ is acyclic in low dimensions,
the cokernel of $d\co C_2(K)\to C_1(K)$ is stably isomorphic to the
augmentation ideal $\fG$. Thus we have a map in $[FC^2(K,X),\fG]$,
which we denote $\nu(\mu)$. Extending Turaev's argument, it is shown
that this coincides with $\nu_{C(K,X),3}(\pd_*(\mu))$. We are now ready
for Bleile's main result.
\begin{theorem}\label{bleilethm}{\rm\cite{bleile}}\qua
Given a triple $(X,G,\mu)$ as above; in particular with $G$ an
(f.p.) group, $X$ incompressible, and a class $\mu\in H_3(K,X;\BZ)$
with $\pd_*(\mu)=[X]$, then if there is an oriented $\mathit{PD}^3$
pair $(Y,X)$ with $\pi_1(Y)=G$ and $i_*[Y]=\mu$, $\nu(\mu)$ is a stable
equivalence.

If also $X$ is incompressible, then conversely if $\nu(\mu)$ is a stable
equivalence, such a $\mathit{PD}^3$ pair exists. Moreover, $(Y,X)$
is unique up to oriented homotopy equivalence.
\end{theorem}
First suppose $(Y,X)$ an oriented $\mathit{PD}^3$ pair with the desired
properties. Then cap product with $[Y]$ gives a chain equivalence
$C^{3-*}(Y,X)\to C_{*}(Y)$, and thus a stable equivalence $\nu[Y]\co
F^2C\to\fG$. Since $i\co Y\to K$ is 2--connected, we may suppose $K$
formed from $Y$ by attaching cells of dimension $\geq 3$, and then
inclusion induces an isomorphism $F^2C(Y,X)\to F^2C(K,X)$. Thus by
naturality, $\nu(i_*[Y])$ also is a stable isomorphism.

Conversely, suppose $\nu(\mu)\in [F^2C(K,X),\fG]$ is a stable equivalence
and $X$ is incompressible. Form $Z$ by attaching 1-- and 2--cells to $X$
and extending the map to $K(G,1)$ until $Z\to K(G,1)$ is 2--connected. Then
$F^2C(K,X)\cong F^2(Z,X)$. The stable equivalence $\nu(\mu)$ thus
arises from an isomorphism $F^2C(Z,X)\oplus P_1\to \fG\oplus P_2$ for
some projective modules $P_1,P_2$, where we may suppose $P_1$ free.
Replacing $Z$ (if necessary) by its bouquet with a number of 2--spheres,
we may suppose $P_1=0$. As before, using an infinite process if necessary,
we may suppose $P_2$ free, of rank $t$, say.

We have an isomorphism of $F^2C(Z,X)$ to $\fG\oplus R^t\subset
R^{t+1}$. Composing with the projection gives a map $C^2(Z,X)\to
R^{t+1}$, which dualises to a map $R^{t+1}\to C_2(Z,X)$ whose image lies
in the kernel of $d_2\co C_2(Z,X)\to C_1(Z,X)$.  Thus the images of the
generators of $R^{t+1}$ define elements of $H_2(Z,X)$.

Since, by the incompressibility hypothesis, $\wtl{X}$ has each
component contractible, the induced map $H_2(Z)\to H_2(Z,X)$ is an
isomorphism. According to our convention, since $\pi_1(Z)=G$, $H_2(Z)$
denotes the homology group of the universal cover of $Z$; thus by
Hurewicz' theorem is isomorphic to $\pi_2(Z)$. Thus each generator of
$R^{t+1}$ defines an element of $\pi_2(Z)$, and we use corresponding maps
$S^2\to Z$ to attach 3--discs to $Z$ to form $Y$. We can then identify
$R^{t+1}=C_3(Y)=C_3(Y,X)$.

It follows as before that $H_3(Y,X;\BZ)\cong \BZ$. Write $[Y]$ for
the generator corresponding to $+1\in\BZ$. Since our construction was
based on the given stable equivalence $\nu(\mu)\in [F^2C(Z),\fG]$,
this equivalence is induced by cap product with $[Y]$, so is equal to
$\nu[Y]$. By naturality, it is also equal to $\nu(i_*[Y])$.  Since
$\wtl{K}$ is contractible and $\wtl{X}$ is homotopy equivalent to a
discrete set, the hypothesis $H_3(K,X)=H_2(K,X)=0$ is satisfied. It
follows by Turaev's Lemma~\ref{turaevlem} that we have $i_*[Y]=\mu$.

In particular, $\pd_*[Y]=\pd_*\mu=[X]$. Hence cap products with $[X]$ and $[Y]$
induce (up to homotopy and sign) a map of exact triangles of chain complexes
$$(C^{2-*}(X)\to C^{3-*}(Y,X)\to C^{3-*}(Y)) \to (C_{*}(X)\to C_{*}(Y)\to
  C_{*}(Y,X)).$$
Since $X$ is a $\mathit{PD}^2$ complex, the map $C^{2-*}(X)\to C_{*}(X)$
is a chain equivalence.  In the induced map of homology sequences
\begin{small}
$$\begin{array}{cccc cccc cccc c}
H^1(X) & \!\!\!\!\!\to\!\!\!\!\! & H^2(Y,X) & \!\!\!\!\!\to\!\!\!\!\! &
  H^2(Y) & \!\!\!\!\!\to\!\!\!\!\! & H^2(X) & \!\!\!\!\!\to\!\!\!\!\! &
  H^3(Y,X) & \!\!\!\!\!\to\!\!\!\!\! & H^3(Y) & \!\!\!\!\!\to\!\!\!\!\! & 0\\
\daw && \daw && \daw && \daw && \daw && \daw && \daw\\
H_1(X) & \!\!\!\!\!\to\!\!\!\!\! & H_1(Y) & \!\!\!\!\!\to\!\!\!\!\! &
  H_1(Y,X) & \!\!\!\!\!\to\!\!\!\!\! & H_0(X) & \!\!\!\!\!\to\!\!\!\!\! &
  H_0(Y) & \!\!\!\!\!\to\!\!\!\!\! & H_0(Y,X) & \!\!\!\!\!\to\!\!\!\!\! & 0
\end{array}$$
\end{small}
we have $H_0(Y)\cong\BZ$ and $H_1(Y)=0$ since this is the homology of
the universal cover of $Y$. Since by construction the map $F^2C(Y,X)\to
C^3(Y,X)$ can be identified with the inclusion $\fG\oplus R^t\subset
R^{t+1}$, we have $H^2(Y,X)=0$ and $H^3(Y,X)\cong\BZ$, so by the
definition of $[Y]$, the induced map $H^3(Y,X)\to H_0(Y)$ is an
isomorphism. Thus all the vertical maps in the diagram are isomorphisms.
As the map of triangles of chain complexes is self-dual, it follows that
we also have isomorphisms for the other half of the sequence, so that
we have chain homotopy equivalences, and so duality holds, as required.

Uniqueness can also be established much as in the former case.

This result enables us to prove the desired splitting theorem.
\begin{theorem}\label{relsplit}
If $(Y,X)$ is a $\mathit{PD}^3$ pair such that $X$ is incompressible and
$G$ is a free product $G'*G''$, then $(Y,X)$ is homotopy equivalent to
the connected sum of $\mathit{PD}^3$ pairs $(Y',X')$ and $(Y'',X'')$,
with $X=X'\,\dot{\cup}\, X''$, $\pi_1(Y')\cong G'$ and $\pi_1(Y'')\cong
G''$. Moreover, this splitting is unique up to homotopy.
\end{theorem}

As before, set $R':=\BZ G'$ etc. For each boundary component $X_r$,
since $H_r=\pi_1(X_r)$ maps injectively to $G$ and is neither free
nor a free product, $H_r$ must be conjugate to a subgroup of $G'$ or
of $G''$. We thus have a natural partition of the components of $X$,
and denote the corresponding partition by $X=X'\,\dot{\cup}\, X''$.

To apply Theorem~\ref{bleilethm} we need a triple $(X',G',\mu')$ with $G'$
an (f.p.)  group, $X'$ incompressible, and $\mu'\in H_3(K',X';\BZ)$ with
$\pd_*(\mu')=[X']$ such that $\nu(\mu')$ is a stable equivalence. We are
given $G'$ and have just constructed $X'$; since $K=K(G,1)$ can be taken
as a bouquet $K'\vee K''$ there is a natural splitting $H_3(K,X;\BZ)\cong
H_3(K',X';\BZ)\oplus H_3(K'',X'';\BZ)$ and we write $\mu=\mu'\oplus
\mu''$. It follows that $\pd_*(\mu')=[X']$, so it remains to discuss
$\nu(\mu')$.

The situation is now essentially the same as before. The class $\nu(\mu')$
is represented by a map, $\phi'\co A'\to B'$ say, of $R'$--modules;
similarly for $\mu''$.  The ring $R$ is free as an $R'$--module, and
since everything splits naturally, $\nu(\mu)$ (which, by hypothesis,
is a stable isomorphism) is represented by $(\phi'\otimes_{R'}R)\oplus
(\phi''\otimes_{R''}R)$. We apply the functor $\otimes_RR'$ and deduce
that $\phi'\oplus((\phi''\otimes_{R'}\BZ)\otimes_{\BZ}R')$ is a stable
isomorphism. We thus need to consider $\phi''\otimes_{R'}\BZ$.

But now, as before, since $\phi''$ is a stable isomorphism of
f.g.\ $R$--modules, $\phi''\otimes_{R'}\BZ$ is a stable isomorphism of
f.g.\ abelian groups. For stable equivalence we can ignore the free
part, and have an isomorphism of the torsion subgroups. Thus, up
to stable equivalence $(\phi''\otimes_{R'}\BZ)\otimes_{\BZ}R'$
is an isomorphism of torsion groups, and the stable isomorphism
$\phi\otimes_RR'$ induces this isomorphism on torsion subgroups,
so the induced map on the torsion-free quotient also is a stable
isomorphism. Hence so is $\phi'$, as required.

\section{$\mathit{PD}^3$ group pairs}\label{secgroup}

We now assume that  $Y$ is a $K(G,1)$; each component $X_r$ of $X$
is a $K(H_r,1)$, say, and $H_r$ injects in $G$. Thus $(G,\{H_r\})$ is
an orientable $\mathit{PD}^3$ pair in the group theoretic sense. An
equivalent notation is to consider a pair $(G,{\bf S})$, where $G$
acts on the set ${\bf S}$: given a collection of subgroups $H_r$ we can
define ${\bf S}$ to be the union of cosets $H\backslash G$, while given a
$G$--set ${\bf S}$ we can choose a set $\{s_i\}$ of orbit representatives
in ${\bf S}$ and define $H_i$ to be the stabiliser of $s_i$. The group $G$
is now torsion free. It is our belief that in this case there is always
a compact 3--manifold $M$ such that $(M,\pd M)$ is homotopy equivalent
to $(K(G,1),\dot{\bigcup}_r K(H_r,1))$. We break this into a series of
lesser problems.

First we need some terminology. A {\it group pair} consists of a group $G$
and a finite list $\{H_r\}$, unordered but perhaps with repetitions, of
subgroups $H_r$.  $(G,\{H_r\})$ is a {\it 2--orbifold pair} if there is
a 2--orbifold $V$ with infinite fundamental group $G$ and the $H_r$ are
the fundamental groups of the boundary components; if $\chi_{orb}(V)<0$
it is a {\it Fuchsian pair}. If there is a normal subgroup $K\lhd G$ with
$K\subseteq \bigcap\{H_r\}$ we have a quotient pair $(G/K,\{H_r/K\})$; if
this is a 2--orbifold pair, we speak of a $K$--by--2--orbifold pair. If
$G$ acts on a tree, with $v$ a vertex, the {\it vertex pair} consists
of the stabiliser $G_v$ of $v$ and the set of stabilisers of the edges
incident to $v$. A pair $(G;G,G)$ is said to be {\it inessential}, a pair
$(G,H)$ with $G$ a torus group and $|G:H|=2$ is {\it weakly inessential}
(this corresponds to the product of $S^1$ with a M\"{o}bius strip);
otherwise we call it essential.

For any field $k$ and group $H$, write $\ol{kH}:=\Hom_k(kH,k)$ for the
set of infinite formal $k$--linear combinations of elements of $H$. If
$H$ is a subgroup of $G$, $G$ acts on $\ol{kG}/(kG\otimes_{kH}\ol{kH})$
and the number of coends $\tl{e}(G,H)$ is the dimension of the submodule
of invariants:
$$\tl{e}(G,H):=\dim_k(\ol{kG}/(kG\otimes_{kH}\ol{kH}))^G.$$
Denote by $\cC_n$ the class of virtually polycyclic groups of Hirsch
rank $n$: thus $\cC_1$ is the class of groups with a subgroup $\BZ$ of
finite index, ie 2--ended groups; $\cC_2$ is the class of groups with
a torus subgroup $\BZ\times \BZ$ of finite index.

We will apply the following major result of Dunwoody and Swenson (there
are similar theorems due to Fujiwara and Papasoglu and to Scott and
Swarup):
\begin{theorem}\label{dunfinal} {\rm\cite{dunswe}}\qua
Let $n\geq 1$,
let $G$ be f.p., and suppose $G$ has no subgroup $H\in \cC_{n-1}$ with
$\tl{e}(G,H)\geq 2$. Then either
\begin{enumerate}
\item[{\rm(i)}] $G\in\cC_{n+1}$,
\item[{\rm(ii)}] $G$ is $\cC_{n-1}$--by--Fuchsian, or
\item[{\rm(iii)}] $G$ has a finite bipartite graph of groups decomposition
such that each black vertex pair is either inessential with $G_v\in
\cC_n$, or a $\cC_{n-1}$--by--Fuchsian pair. (Thus every edge group is in
$\cC_n$.)
\end{enumerate}
Moreover,
\begin{enumerate}
\item[{\rm(a)}] For every splitting $G=A*_CB$ or $G=A*_C$ of $G$ over
a $\cC_n$ group $C$, each white vertex group is conjugate into $A$ or $B$.
\item[{\rm(b)}] For every $\cC_{n-1}$--by--Fuchsian pair which is a
vertex pair in some splitting of $G$, its vertex group is conjugate into
a black vertex group.
\item[{\rm(c)}] Every subgroup $J\in\cC_n$ of $G$ with  $\tl{e}(G,J)\geq
2$ is conjugate into a vertex group and has a subgroup of finite index
conjugate into a black vertex group.
\end{enumerate}
\end{theorem}

\begin{theorem}\label{splitPD3}
If $G$ is an f.p.\ orientable $\mathit{PD}^3$ group, one of the following holds.
\begin{enumerate}
\item[{\rm(i)}] $G\in\cC_3$.
\item[{\rm(ii)}] $G$ is a $\BZ$--by--Fuchsian group.
\item[{\rm(iii)}] $G$ has a finite bipartite graph of groups decomposition such
that each edge group is a torus group; for each black vertex, either
the vertex group is a torus group and the vertex pair inessential,
or the vertex pair is a $\BZ$--by--Fuchsian pair.
\end{enumerate}
Moreover,
\begin{enumerate}
\item[{\rm(a)}] For every splitting $G=A*_CB$ or $G=A*_C$ of $G$ over
a torus group $C$, each white vertex group is conjugate into $A$ or $B$.
\item[{\rm(b)}] For every $\BZ$--by--Fuchsian pair which is a vertex
pair in some splitting of $G$, its vertex group is conjugate into a
black vertex group.
\item[{\rm(c)}] Every torus subgroup $J$ of $G$ is conjugate into a
vertex group and has a subgroup of finite index conjugate into a black
vertex group.
\end{enumerate}
\end{theorem}
\begin{proof} Since $G$ is a $\mathit{PD}^3$ group, we have $\cd G=3$;
in particular, $G$ is torsion free (so any $\cC_1$ subgroup is
isomorphic to $\BZ$). It also follows that for a subgroup $H\cong\BZ$,
$\tl{e}(G,H)=1$. Thus the hypotheses of Theorem~\ref{dunfinal} are
satisfied, taking $n=2$. Since a torsion-free $\cC_2$ group is either
a torus group or a Klein bottle group, and as $G$ is orientable it
cannot split over a Klein bottle group, and since for any torus subgroup
$J\subset G$, $\tl{e}(G,J)=2$, the result follows.
\end{proof}
\begin{addendum}
In case (iii), we also have
\begin{enumerate}
\item[{\rm(d)}] every torus subgroup $J$ of $G$ is conjugate into a
black vertex group,
\item[{\rm(e)}] each essential white vertex pair is atoroidal.
\end{enumerate}
\end{addendum}
\begin{proof}
To prove (d), we may suppose $J$ contained in a white vertex group. Thus
(d) will follow from (e).

Denote the white vertex pair in question by $(G',{\bf S})$. By (c),
any torus subgroup $J\subset G'$ has a subgroup $J'$ of finite index
conjugate into a black vertex group. Thus $J'$ is conjugate into two
distinct vertex groups of the decomposition. It follows from properties
of subgroups of split groups that $J'$ is contained in one of the edge
groups $H'\in {\bf S}$.

We may suppose without loss of generality that $J$ is a maximal torus
subgroup of $G'$, and that $J'$ is a maximal subgroup which is conjugate
into $H'$. Visualise $(G',{\bf S})$ as a manifold pair $(M,\pd M)$;
consider the covering space $\tl{M}$ corresponding to the subgroup
$J'$. The boundary components over $H'$ correspond to double cosets of
$H'$ and $J'$ in $G'$; to the double coset $J'gH'$ corresponds something
with fundamental group $g^{-1}J'g\cap H'$. If $g\in J$, this is just
the torus group $J'$. If the cosets $J'g,J'g'$ (with $g,g'\in J$) are
distinct, so are $J'gH'$ and $J'g'H'$, for otherwise there exists $h\in
H'$ with $J'gh=J'g'$ and so $h\in J$; by maximality of $J'$, $h\in J'$,
so $J'g=J'g'$. Thus in the cover $\tl{M}$ --- itself homotopy equivalent
to a torus --- there are at least $|J:J'|$ boundary components which
are tori. Hence $H_3(\tl{M},\pd\tl{M})$ has rank at least $|J:J'|-1$. By
duality, so has $H^0_c(\tl{M})$. But this group vanishes if the covering
is infinite and $\tl{M}$ non-compact; if the cover is finite, the group
has rank 1. Hence if $|J:J'|>1$, the index is 2; the covering is finite,
so $J'$ has finite index in $G'$; as $G'$ is torsion free, it is a torus
group, so $G'=J$. Now $H'=J'$ is a torus group of index 2 in $G'$,
thus $(G',H')$ is weakly inessential.  Otherwise $J=J'$ is contained
in a boundary group. Since this holds for any torus subgroup of $G'$,
the vertex pair is atoroidal.
\end{proof}
In the case when $Y=K(G,1)$ is a 3--manifold, one can show that the
essential white vertex pairs are not only atoroidal but acylindrical. We
expect this to be the case here also.
\begin{conjecture}\label{acylind}
Let $(G;\{H_r\})$ be an essential white vertex pair in the above
decomposition, and suppose that $g\in G, s_i\in H_r$ and $s_j\in H_j$
satisfy $g^{-1}s_ig=s_j$. Then either $s_i=1$ (and hence $s_j=1$) or $i=j$
and $g^{-1}H_rg=H_r$.\end{conjecture}

Peter Scott indicated to me a proof in the geometric case. Trying to
push it through in the $\mathit{PD}$ case leads me to the following.

Suppose $s$ belongs to two different edge groups. Since each edge group
is a torus group, we can choose elements $t_1$, $t_2$ in them such that
$T_i=\la s,t_i\ra$ is a subgroup of finite index in the edge group. Define
$t_1t_2t_3=1$. Then $s$ also commutes with $t_3$, so the group $T_3$ they
generate is either cyclic or is contained in a further edge group. Suppose
$H:=\la s,t_1,t_2,t_3\ra$ is the direct product of $\la s\ra$ and a free
group on the $t_i$. Then $(H;T_1,T_2,T_3)$ is a $\mathit{PD}^3$ pair (the
`product of a circle and a pair of pants'). The covering map from $K(H,1)$
to $K(G,1)$ gives a map of Poincar\'e pairs, of finite degree on each
boundary component, hence of finite degree. Thus $H$ has finite index
in $G$. But $H$ is not atoroidal (there are too many torus subgroups),
hence nor is $G$.

We now wish to extend the above results to the case $X\neq\emptyset$. If
$G$ satisfies the maximal condition on centraliser subgroups (denoted
Max--c), a theorem of Kropholler in the case $n=3$ gives:
\begin{theorem} {\rm\cite{krop}}\qua
Let $(G,{\bf S})$ be a $\mathit{PD}^3$--pair such that $G$ satisfies
Max--c but is not in $\cC_3$. Then there is a unique reduced $G$--tree
$Y$, adapted to ${\bf S}$, such that $G\backslash Y$ is finite, each
edge group is in $\cC_2$ (hence is a torus group), each vertex pair is
either of Seifert type or atoroidal; and every torus subgroup of $G$
fixes a vertex of $Y$.
\end{theorem}

We expect that the hypothesis Max--c is unnecessary, more
precisely: \begin{conjecture}\label{pairsplit} There is an analogue
of Theorem~\ref{dunfinal} for group pairs.  \end{conjecture} It would
suffice for us to have an extension of Theorem~\ref{splitPD3}. One line
of attack is to apply Theorem~\ref{splitPD3} to the double $D(G,{\bf
S})$. Since the resulting splitting is unique (this follows from a theorem
of \cite{foruniq}, as apart from inessential components, the tree is
`strongly slide-free'), it is compatible with the involution $\gs$
interchanging the two copies of $G$. It would be necessary to investigate
how $\gs$ restricts to the vertex pairs of the splitting.

I abstain from further discussion, as it seems likely that this
question can be resolved by the methods used for the proof of
Theorem~\ref{dunfinal}.

\section{Geometrisation}\label{secgeom}

The hope is to find some sort of decomposition of any $\mathit{PD}^3$
pair into pieces which are geometric in some homotopy-theoretic sense. We
refer to \cite{scott} for an account of the Thurston geometries.

In sections \ref{secsphere} and \ref{secdisc} we discussed decomposition
of an arbitrary pair into pieces with $G$ not a free product: while this
does not hold in all cases, it does if $G$ is torsion free. Suppose now
$G$ not a free product. If $e(G)=2$ we have $S^2\times S^1$ and $D^2\times
S^1$. If $e(G)=0$, we have complexes $Y$ finitely covered (up to homotopy)
by $S^3$. Each of these may be regarded as geometric. The group $G$ has
periodic cohomology with period 4. Such groups have been classified, and
each choice of a generator $g\in H^4(G;\BZ)$ determines a $\mathit{PD}^3$
complex, unique up to homotopy equivalence.

In the case $e(G)=1$, we gave a further decomposition in
Theorem~\ref{splitPD3}. In cases (i) and (ii), $K(G,1)$ is homotopy
equivalent to a 3--manifold with geometric structure; of type $\BR^3$,
$\mathit{Nil}$ or $\mathit{Sol}$ for (i) and $\BH^2\times \BR$ or
$\wtl{SL_2(\BR)}$ for (ii). In case (iii) the black vertex pairs
correspond to Seifert manifolds with boundary, which have geometric
structures of Seifert type.  This identification is essentially due to
the theorem that convergence groups are Fuchsian. If Hope~\ref{pairsplit}
is true, corresponding results also hold in the case of pairs.

Much harder is the case of the essential white vertex pairs. It was
conjectured by Kropholler that if $(G,{\bf S})$ is atoroidal, $G$
is isomorphic to a discrete subgroup of $PSL_2(\BC)$ and ${\bf S}$
to the collection of peripheral subgroups. We attempt a more precise
formulation as two separate questions.
\begin{conjecture}\label{hypcjc}
If the $\mathit{PD}^3$ pair $(G,{\bf S})$ is atoroidal and acylindrical,
then $G$ is hyperbolic relative to $\mathbf{S}$.
\end{conjecture}
Misha Kapovich and Bruce Kleiner have recently shown that if $G$
is a $\mathit{PD}^3$ group which acts discretely, isometrically and
cocompactly on a CAT(0) space, then $G$ is either a hyperbolic group,
a Seifert manifold group, or splits over a virtually abelian subgroup.
\begin{conjecture}\label{geomcjc}
If the $\mathit{PD}^3$ pair $(G,{\bf S})$ is atoroidal and acylindrical,
and $G$ is hyperbolic relative to $\mathbf{S}$, then $G$ is isomorphic
to a discrete subgroup of $PSL_2(\BC)$ and ${\bf S}$ to the collection
of peripheral subgroups.
\end{conjecture}

Cannon and co-workers have written several papers developing an approach
to this problem in the case ${\bf S}=\emptyset$; I confine myself to
citing \cite{can} and \cite{canswe}, and the outline in \cite{kapben}. It
is known that a hyperbolic group is a $\mathit{PD}^3$ group if and
only if $\pd G$ has the \v{C}ech homology of a 2--sphere; and that
then $\pd G$ is homeomorphic to $S^2$. If there is a homeomorphism
preserving the quasi-conformal structure, or equivalently if $G$
is itself quasi-isometric to hyperbolic space $\BH^3$, a theorem of
Sullivan applies to give the desired result.

If ${\bf S}\neq\emptyset$, then $\cd G=2$. Hence $\dim\pd G=1$. It
follows from duality that $G$ cannot split over a finite group, so $\pd
G$ is connected; it follows from the acylindricity hypothesis that $\pd
G$ has no local cut points. It now follows from a theorem of Kapovich
and Kleiner \cite{KK} that $\pd G$ is homeomorphic to either $S^1$,
a Sierpinski gasket or a Menger curve.

Consider the double $D(G,{\bf S})$. Either $G$ is a surface group or
$D(G,{\bf S})$ is hyperbolic, with $G$ a quasiconvex subgroup. Hence
$\pd G$ embeds in $\pd D(G,{\bf S})$, which is homeomorphic to $S^2$,
so is planar. Thus $\pd G$ cannot be a Menger curve, so is either $S^1$
or the Sierpinski gasket, and if $\pd G=S^1$, $G$ is Fuchsian. The
Sierpinski gasket is obtained from a 2--disc by removing the interiors
of a sequence of disjoint discs, so has a natural boundary which is a
union of circles $S_r$. In this case Kapovich and Kleiner further show:
\begin{theorem}\label{gasket}
The circles $S_r$ fall into finitely many orbits under $G$. The stabiliser
of $S_r$ is a virtually Fuchsian group $H_r$, quasi-convex in $G$,
which acts on $S_r$ as a uniform convergence group.

Let ${\bf H}$ consist of one representative from each conjugacy class
of subgroups $H_r$. Then the double $D(G,{\bf H})$ is hyperbolic, $G$
is a quasiconvex subgroup of it, $\pd D(G,{\bf H})$ is homeomorphic
to $S^2$, and $D(G,{\bf H})$ is a $\mathit{PD}^3$ group (over $\BQ$;
if torsion-free, over $\BZ$).
\end{theorem}
Thus if Hope~\ref{geomcjc} can be established when $X=\emptyset$, the
general case seems to follow.

\Addresses\recd

\end{document}